\documentclass[12pt]{article}
\usepackage{latexsym,amsfonts,amsmath,amssymb}
\newtheorem{theorem}{Theorem}
\def\Theorem. #1\par{\begin{theorem}#1\end{theorem}}
\newtheorem{corollary}[theorem]{Corollary}
\def\Corollary. #1\par{\begin{corollary}#1\end{corollary}}
\newtheorem{lemma}[theorem]{Lemma}
\def\Lemma. #1\par{\begin{lemma}#1\end{lemma}}
\newtheorem{question}[theorem]{Question}
\def\Question. #1\par{\begin{question}#1\end{question}}
\newtheorem{observation}[theorem]{Observation}

\newtheorem{claim}[theorem]{Claim}
\def\Claim. #1\par{\begin{claim}#1\end{claim}}
\newtheorem{conjecture}[theorem]{Conjecture}
\def\Conjecture. #1\par{\begin{conjecture}#1\end{conjecture}}
\newtheorem{definition}[theorem]{Definition}
\def\Definition. #1\par{\begin{definition}#1\end{definition}}
\newtheorem{remark}[theorem]{Remark}
\def\Remark. #1\par{\begin{remark}#1\end{remark}}
\newcommand{\F}{{\mathbb F}}
\def\P{{\mathbb P}}
\newcommand{\Q}{{\mathbb Q}}
\newcommand{\Qdot}{\dot\Q}
\newcommand{\Ftail}{\F_{\!\scriptscriptstyle\rm tail}}
\newcommand{\Fotail}{\F^0_{\!\scriptscriptstyle\rm tail}}

\newcommand{\fotail}{f^0_{\!\scriptscriptstyle\rm tail}}
\newcommand{\Gtail}{G_{\!\scriptscriptstyle\rm tail}}
\newcommand{\Gotail}{G^0_{\!\scriptscriptstyle\rm tail}}
\newcommand{\Ptail}{\P_{\!\scriptscriptstyle\rm tail}}
\newcommand{\Potail}{\P^0_{\!\scriptscriptstyle\rm tail}}
\newfont{\msam}{msam10 at 12pt}
\newcommand{\from}{\mathbin{\vbox{\baselineskip=3pt\lineskiplimit=0pt
                         \hbox{.}\hbox{.}\hbox{.}}}}
\newcommand{\of}{\subseteq}
\newcommand{\set}[1]{\left\{\,{#1}\,\right\}}
\newcommand{\compose}{\circ}
\newcommand{\elesub}{\prec}
\newcommand{\muchgt}{\gg}

\newcommand{\inverse}{{-1}}
\newcommand{\dom}{\mathop{\rm dom}}
\newcommand{\ran}{\mathop{\rm ran}}

\newcommand{\Add}{\mathop{\rm Add}}
\newcommand{\image}{\mathbin{\hbox{\tt\char'42}}}
\newcommand{\plus}{{+}}
\newcommand{\plusplus}{{{+}{+}}}
\newcommand{\restrict}{\mathbin{\hbox{\msam\char'26}}}
\newcommand{\satisfies}{\models}
\newcommand{\forces}{\Vdash}

\newcommand{\cross}{\times}
\newcommand{\union}{\cup}

\newcommand{\intersect}{\cap}

\newcommand{\ltkappa}{{{<}\kappa}}

\newcommand{\leqgamma}{{{\leq}\gamma}}
\newcommand{\ltomega}{{{<}\omega}}
\newcommand{\leqtheta}{{{\leq}\theta}}
\newcommand{\leqbeta}{{{\leq}\beta}}

\newcommand{\leqdelta}{{{\leq}\delta}}
\newcommand{\card}[1]{\left|#1\right|}
\def\boolval#1{\mathopen{\lbrack\!\lbrack}\,#1\,\mathclose{\rbrack\!
        \rbrack}}
\newcommand{\st}{\mid}

\newcommand{\iso}{\cong}
\def\<#1>{\langle\,#1\,\rangle}

\newcommand{\QEDbox}{\fbox{}}
\newcommand{\factordiagramup}[6]{$$\begin{array}{ccc}
#1&\raise3pt\vbox{\hbox to60pt{\hfill$\scriptstyle #2$\hfill}\vskip-6pt\hbox{$\vector(4,0){60}$}}&#3\\
\vbox to30pt{}&\raise22pt\vtop{\hbox{$\vector(4,-3){60}$}\vskip-22pt\hbox to60pt{\hfill$\scriptstyle #4\qquad$\hfill}}
     &\ \ \lower22pt\hbox{$\vector(0,3){45}$}\ {\scriptstyle #5}\\
\vbox to15pt{}&&#6\\
\end{array}$$}
\newcommand{\factordiagram}[6]{$$\begin{array}{ccc}
#1&&\\
\ \ \raise22pt\hbox{$\vector(0,-3){45}$}\ {\scriptstyle #2}
&\raise22pt\hbox{$\vector(2,-1){90}$}\raise5pt\llap{$\scriptstyle#3$\qquad\quad}&\vbox to25pt{}\\
#4&\raise3pt\vbox{\hbox to90pt{\hfill$\scriptstyle #5$\hfill}\vskip-6pt\hbox{$\vector(4,0){90}$}}&#6\\
\end{array}$$}
\newcommand{\df}{\it} 
\newcommand{\beginnamedtheorem}[1]{\newtheorem{#1}[theorem]{#1}\begin{#1}}
\newenvironment{proof}{\noindent{\bf Proof: }}{\kern2pt\QEDbox\par\bigskip}
\def\Proof{\begin{proof}}
\newcommand{\QED}{\end{proof}}
\begin{document}
\author{Joel David Hamkins\thanks{My research has been supported in part by grants from the PSC-CUNY Research Foundation and the NSF. I would particularly like to thank Philip Welch for suggesting this line of research to me and for his helpful discussions concerning it during my recent return to Japan in July, 1999. And I would like also to thank Kobe University for their generous support during that trip.}\\
{\normalsize\sc City University of New York\thanks{Specifically, the College of Staten Island of CUNY and the CUNY Graduate Center.}}\\
{\footnotesize http://www.math.csi.cuny.edu/$\sim$hamkins}\\}
\title{Unfoldable cardinals and the GCH}
\maketitle

\abstract{\noindent Unfoldable cardinals are preserved by fast function forcing and the Laver-like preparations that fast functions support. These iterations show, by set-forcing over any model of ZFC, that any given unfoldable cardinal $\kappa$ can be made
indestructible by the forcing to add any number of Cohen subsets to $\kappa$.}

\section{Unfoldable Cardinals}

In introducing unfoldable cardinals last year, Andres Villaveces \cite{Villaveces98} ingeniously extended the notion of weak compactness to a larger context, thereby producing a large cardinal notion, unfoldability, with some of the feel and flavor of weak compactness but with a greater consistency strength. In this paper I will show that the embeddings associated with these unfoldable cardinals are amenable to some of the same lifting techniques that apply to weakly compact embeddings, augmented with some methods from the strong cardinal context. Using these techniques, I will show by set-forcing that any unfoldable cardinal $\kappa$ can be made indestructible by the forcing to add any number of Cohen subsets to $\kappa$.

Villaveces defines that a cardinal $\kappa$ is {\df $\theta$-unfoldable} when for every transitive model $M$ of size $\kappa$ with $\kappa\in M$ there is an elementary embedding $j:M\to N$ with critical point $\kappa$ and $j(\kappa)\geq\theta$.\footnote{Actually, Villaveces defines that $\kappa$ is $\theta$-unfoldable if and only if it is inaccessible and for every $S\of\kappa$ there is an $\hat S$ and a transitive $N$ of height at least $\theta$ such that $\<V_\kappa,{\in},S>\elesub\<N,{\in},\hat S>$. He then proves this definition equivalent to the embedding characterization, which I prefer to take as the basic notion.} The cardinal $\kappa$ is fully {\df unfoldable} when it is $\theta$-unfoldable for every ordinal $\theta$. A lower bound on the consistency strength of unfoldability is provided by the observation that $\kappa$ is weakly compact if and only if $\kappa$ is $\kappa$-unfoldable. For an easy upper bound, one can see by simply iterating a normal measure that every measurable cardinal is unfoldable; indeed, Villaveces proves that every Ramsey cardinal is unfoldable. 

Villaveces is interested in the question of when an unfoldable cardinal is preserved by certain forcing iterations, and in particular when the GCH can be forced to fail at a given unfoldable cardinal. He defines $\kappa$ to be {\df strongly unfoldable} when for every $\theta$ and every transitive set $M\satisfies ZF^-$ of size $\kappa$ with $\kappa\in M$ and $M^\ltkappa\of M$ there is an embedding $j:M\to N$ with critical point $\kappa$ such that $j(\kappa)\geq\theta$ and $V_\theta\of N$, and then conjectures:

\Conjecture. {\rm (Villaveces \cite{Villaveces98})} If the existence of a strongly unfoldable cardinal is consistent with ZFC, then so is the existence of an unfoldable cardinal $\kappa$ with $2^\kappa>\kappa^\plus$. 

Since any unfoldable cardinal is unfoldable in $L$, and any unfoldable cardinal in $L$ is strongly unfoldable there (since $V_\theta^L\of L_{\aleph_\theta}$), the existence of an unfoldable cardinal is equiconsistent with the existence of a strongly unfoldable cardinal. Thus, the notion of strong unfoldability in the conjecture can be completely avoided. Villaveces is really conjecturing that unlike the situation for measurable cardinals, having the GCH fail at an unfoldable cardinal does not increase consistency strength. This conjecture is a consequence of the following theorem, the main result of this paper.

\beginnamedtheorem{Main Theorem} If $\kappa$ is unfoldable in $V$, then there is a forcing extension, by forcing of size $\kappa$, in which the unfoldability of $\kappa$ becomes indestructible by $\Add(\kappa,\theta)$ for any $\theta$.
\end{Main Theorem}

Villaveces and Leshem \cite{Villaveces&Leshem} have also announced a proof of the conjecture, although not of the theorem I have just claimed here. Given Villaveces' embedding characterization of unfoldability, though, their proof is not the kind of lifting argument that one might initially expect; indeed, their proof employs proper class forcing over $L$ and does not use embeddings at all. In my proof, by set forcing over any model with an unfoldable cardinal, I try to remain close to the embeddings. Indeed, borrowing some techniques from the strong cardinal context in order to treat unfoldability embeddings much like weak compactness embeddings, I simply lift the unfoldability embeddings to the forcing extension. And by doing so, I hope to have provided a general lifting technique for unfoldable cardinals, a technique that I expect will be easily adapted to many other purposes.

Before proceeding further, I should mention a claim in Villaveces and Leshem's \cite{Villaveces&Leshem} that could prove potentially troublesome, namely, their Theorem 2.1, which asserts that set forcing cannot achieve what I have just claimed of it in the Main Theorem. Specifically, they claim that over $L$, any set forcing making $2^\kappa>\kappa^\plus$ will destroy the unfoldability of $\kappa$. My Main Theorem asserts that by set forcing over any model with an unfoldable cardinal --- including $L$, and any unfoldable cardinal is unfoldable in $L$ --- one can force $2^\kappa$ to be as large as desired while preserving the unfoldability of $\kappa$. So there is some conflict between our claims. Simply put, either one of us is wrong, or there are no unfoldable cardinals. I believe that the error of Villaveces and Leshem's Theorem 2.1 is their claim that $M\of L_\rho[G]$ in their argument. The corresponding $M$ provided by the proof of my Main Theorem, and I will point this out again just after the proof, is defined in $L[G]$ but not in $L_\rho[G]$. Thus, my proof provides a counterexample to this specific step in their argument.

Before beginning the proof of the theorem, let me make a few observations about unfoldable cardinals in order to simplify the arguments to come. I have said that a cardinal $\kappa$ is $\theta$-unfoldable when every transtive set $M$ of size $\kappa$ with $\kappa\in M$ has an embedding $j:M\to N$ with critical point $\kappa$ and $j(\kappa)\geq\theta$. The first observation is that if $M$ contains another structure $M'$, then $j\restrict M'$ provides a $\theta$-unfoldability embedding from $M'$ into $j(M')$. Consequently, to verify the unfoldability of $\kappa$ it suffices to consider only the members of some unbounded family of such $M$. Specifically, if for some enormous $\eta\muchgt\kappa$ we have an elementary substructure $X\elesub V_\eta$ of size $\kappa$ with $\kappa+1\of X$ and $X^\ltkappa\of X$ (and any transitive structure $M'$ of size $\kappa$ can be placed into such an $X$), then let $M$ be the Mostowski collapse of $X$. It follows that $M$ has size $\kappa$, $M^\ltkappa\of M$ and (by choosing $\eta$ appropriately) $M\satisfies ZFC^-$, a suitably large fragment of ZFC. In this paper, I will define such $M$ as the {\df nice} models of $ZFC^-$ (of size $\kappa$). If $X$ did contain as an element some transitive structure $M'$ of size $\kappa$, then since $M'$ would be fixed by the collapse, it follows that $M'\in M$. Thus, when verifying the unfoldability of $\kappa$, it suffices only to have unfoldability embeddings for the nice models.

My second observation concerns a canonical form for the $\theta$-unfoldability embeddings $j:M\to N$ for such nice models $M$. Given such an embedding, let $X=\set{j(g)(s) \st g\in M\And s\in\theta^\ltomega}$. By verifying the Tarski-Vaught criterion, it is easy to see that $X$ is an elementary substructure of $N$; indeed, it is what I have called elsewhere (see \cite{Seeds}) the {\df seed hull} of $\theta^\ltomega$, the Skolem hull in $N$ of all seeds $s\in\theta^\ltomega$ with $\ran(j)$. If $\pi:X\iso N_0$ is the Mostowski collapse of $X$, then the embedding $j_0=\pi\compose j$ provides an elementary embedding $j_0:M\to N_0$ and the following commutative diagram, where $h=\pi^\inverse$. 

\factordiagram{M}{j_0}{j}{N_0}{h}{N}

\noindent Since $\theta\of X$ if follows that $\pi$ is the identity on $\theta$, and so $j_0(\kappa)=\pi(j(\kappa))$ remains at least as large as $\theta$. Also, because of the form of $X$ it follows easily from this that $N_0=\set{j_0(g)(s) \st g\in M\And s\in\theta^\ltomega}$. Thus, $j_0:M\to N_0$ is a $\theta$-unfoldability embedding whose target space $N_0$ is the hull of $\theta^\ltomega$ and $\ran(j_0)$. 

In particular, if $\kappa$ is $(\theta+1)$-unfoldable, and this just means $j(\kappa)>\theta$, then for every nice model $M$ of $ZFC^-$, there is an embedding $j:M\to N$ with $j(\kappa)>\theta$ and $N=\set{j(g)(s) \st g\in M\And s\in(\theta+1)^\ltomega}$. I will refer to such embeddings as the {\df canonical} $(\theta+1)$-unfoldability embeddings.

\section{Unfoldability after fast function forcing}\label{FastFunctionSection}

I would now like to show that the unfoldability of an unfoldable cardinal $\kappa$ is preserved by Woodin's fast function forcing, the notion of forcing that adds a generic partial function $f\from\kappa\to\kappa$ having many of the features of a Laver function in the supercompact cardinal context. This notion of forcing has proved useful in other large cardinal contexts, particularly that of strongly compact cardinals. In \cite{LotteryPrep}, I showed that weakly compact cardinals, measurable cardinals, strong cardinals, strongly compact cardinals and supercompact cardinals are all preserved by fast function forcing, and furthermore in each of these large cardinal contexts the fast function leads to a generalized notion of Laver function appropriate for that context. In that paper, I used these functions to provide a new general kind of Laver preparation, what I called the lottery preparation, for each of these kinds of cardinals, making them indestructible to varying extents. So while I have included all the necessary details for fast function forcing in this paper, I refer the reader to \cite{LotteryPrep} for a more gradual, general introduction. The basic feature of fast function forcing that makes it so useful is that in each of the different large cardinal contexts, one can have embeddings $j:V\to M$ (or, here, $j:M\to N$) for which $j(f)(\kappa)$ is almost arbitrarily specified. In this paper, I am pleased to add unfoldable cardinals to the list, including them within the realm of fast function forcing's efficacy. 

Conditions in the fast function forcing poset $\F$ are the partial functions $p\from\kappa\to\kappa$ of size less than $\kappa$ such that if $\gamma$ is in $\dom(p)$, then $\gamma$ is inaccessible, $p\image\gamma\of\gamma$ and $\card{p\restrict \gamma}<\gamma$. The conditions are ordered by inclusion. Forcing with $\F$ adds a partial function $f\from\kappa\to\kappa$, the {\df fast function}, with the property that every element in $\dom(f)$ is an inaccessible closure point of $f$. Below the condition $p=\set{\<\beta,\alpha>}$, by splitting each condition into small and large parts, the forcing $\F\restrict p$ factors as $\F_\beta\cross\F_{\gamma,\kappa}$, where $\F_\beta$ consists of those conditions of size less than $\beta$ with domain contained in $\beta$, and $\F_{\gamma,\kappa}$ consists of those conditions with domain in $(\gamma,\kappa)$ where $\gamma=\max\set{\kappa,\alpha}$. It is relatively easy to see, by simply taking the union of conditions, that $\F_{\gamma,\kappa}$ is $\leqgamma$-closed. 

\Theorem. Any unfoldable cardinal $\kappa$ is preserved by fast function forcing. Specifically, every canonical $(\theta+1)$-unfoldability embedding $j:M\to N$ in the ground model lifts to a $(\theta+1)$-unfoldability embedding $j:M[f]\to N[j(f)]$ in the extension $V[f]$. Furthermore, if $\alpha$ is any ordinal below $j(\kappa)$, then the lift can be chosen in such a way that $j(f)(\kappa)=\alpha$.\label{FastFunction}

\Proof First, let me quickly point out that in order to verify unfoldability in the extension it suffices merely to lift the ground model embeddings. This is because if $M'$ is any structure of size $\kappa$ in $V[f]$, then it has a name $\dot M'$ of size $\kappa$ in $V$, and this name can be placed into the domain of a ground model canonical unfoldability embedding $j:M\to N$. If this ground model embedding lifts to $j:M[f]\to N[j(f)]$, then $M'=\dot M'_f$ is included in $M[f]$, the domain of an unfoldability embedding in $V[f]$, and this is sufficient, since as I remarked earlier the restriction $j\restrict M':M'\to j(M')$ then provides an unfoldability embedding for $M'$. 

So fix any $\theta$ and any canonical $(\theta+1)$-unfoldability embedding $j:M\to N$ in $V$, and choose any $\alpha<j(\kappa)$. I will show how to lift the embedding to the extension. For simplicity, let me initially assume that $\theta\leq\alpha$; I will consider the general case subsequently. Below the condition $p=\set{\<\kappa,\alpha>}$, a partial function with only one element in its domain, we may split the conditions of $j(\F)$ into their small part, with domain below $\kappa$, and their tail part, with domain above $\alpha$. Since the collection of such small parts is precisely $\F$ again, this splitting procedure shows that $j(\F)\restrict p$ is isomorphic to $\F\cross\Ftail$. In order to lift the embedding, therefore, my strategy will first be to find a generic for $\Ftail$. 

To do so, I will employ a factor technique of W. Hugh Woodin's from the strong cardinal context. Specifically, let $X=\set{ j(g)(\kappa,\theta,\alpha) \st g\in M }$, the seed hull of $\<\kappa,\theta,\alpha>$ in $N$. This is an elementary substructure of $N$. Let $\pi:X\iso N_0$ be the Mostowski collapse, and factor the embedding to produce $j_0:M\to N_0$ and $h:N_0\to N$, where $j_0=\pi\compose j$, $h=\pi^\inverse$ and $j=h\compose j_0$. 

\factordiagram{M}{j_0}{j}{N_0}{h}{N}

Because there are only $\kappa$ many functions $g$ in $M$, it follows that $X$, and hence also $N_0$, have size $\kappa$. Let $p_0=\pi(p)=\set{\<\kappa,\alpha_0>}$, where $h(\alpha_0)=\alpha$. Let 
$\theta_0=\pi(\theta)$ so that $h(\theta_0)=\theta$. Below $p_0$, the 
forcing $j_0(\F)$ factors as $\F\cross\Fotail$, where $\Fotail$ is ${\leq}\theta_0$-closed in $N_0$. 

I claim that $N_0^\ltkappa\of N_0$ in $V$. To see why, suppose that $\<a_\alpha\st\alpha<\beta>$ is in $N_0^\ltkappa$. Each $a_\alpha$, being in $N_0$, must have the form $a_\alpha=j_0(g_\alpha)(\kappa,\theta_0,\alpha_0)$ for some function $g_\alpha\in M$. Since $M^\ltkappa\of M$, the sequence $\vec g=\<g_\alpha\st\alpha<\beta>$ is also in $M$, so 
the sequence $j_0(\vec g)=\<j_0(g_\alpha)\st\alpha<\beta>$ is in $N_0$. Thus, by simply evaluating these functions at $(\kappa,\theta_0,\alpha_0)$, the sequence $\<a_\alpha\st\alpha<\beta>$ is also in $N_0$, and I have established my claim.

Continuing now with the main line of argument, since $N_0^\ltkappa\of N_0$, we may line up the $\kappa$-many dense subsets of $\Fotail$ in $N_0$ into a $\kappa$-sequence in $V$ and diagonalize 
against these dense sets to produce in $V$ an $N_0$-generic $\fotail$. Specifically, we recursively construct a descending sequence of conditions which gradually meets more and more of the dense sets; at limit stages, the sequence constructed so far is in $N_0$, since $N_0^\ltkappa\of N_0$, and there is a condition below it by the closure of the forcing $\Fotail$ in $N_0$. So the constuction continues up to $\kappa$. If $\fotail$ is the filter generated by this descending sequence, then, it is constructed in $V$ but it meets all the required dense sets to be $N_0$-generic for $\Fotail$. Since $f$ is $V$-generic it is also $N_0[\fotail]$-generic, and so the function $f\union p_0\union \fotail$ is 
$N_0$-generic for $j_0(\F)$. Thus, $j_0$ lifts to $j_0:M[f]\to N_0[j_0(f)]$ with 
$j_0(f)=f\union p_0\union \fotail$.

Moving over to $N$, now, I will argue that the filter generated by $h\image \fotail$ is $N$-generic for $\Ftail$.
Every element of $N$ has the form $j(g)(s)$ for 
some function $g$ in $M$ and some $s$ in $\theta^\ltomega$. Thus, since 
$j(g)=h(j_0(g))$, every element of $N$ has the form $h(r)(s)$ for some 
function $r$ in $N_0$. And we may assume that $r:\theta_0^{|s|}\to N$, 
since this is enough to ensure that $s$ is in $\dom(h(r))$. Thus, if $D$ is a 
dense open subset of $\Ftail$ in $N$, then $D$ must be $h(r)(s)$ for some 
such $r$. We may assume without loss of generality that $r(\sigma)$ is a 
dense open subset of $\Fotail$ for every $\sigma$ in $\theta_0^{|s|}$. But 
$\Fotail$ is actually ${\leq}\theta_0$-closed in $N_0$, so we may intersect 
all these $r(\sigma)$ to obtain a single dense set $D^*$ in $N_0$ which is 
contained in $r(\sigma)$ for every $\sigma$. It follows that $h(D^*)$ is 
contained in $D$. Also, since $D^*$ meets $\fotail$ by genericity, it follows that $h(D^*)$ meets 
$h\image \fotail$. Thus, $D$ meets $h\image\fotail$, as desired. 

Since $f$ is $V$-generic, it is also $N[(h\image\fotail)]$-generic, and so 
$f\union p\union(h\image\fotail)$ is $N$-generic for $j(\F)$. Thus, we may lift the embedding 
$j$ to $j:M[f]\to N[j(f)]$ where $j(f)=f\union p\union(h\image\fotail)$, resulting in the following diagram.

\factordiagram{M[f]}{j_0}{j}{N_0[j_0(f)]}{h}{N[j(f)]}

\noindent In particular, $j(f)(\kappa)=p(\kappa)=\alpha$, as desired. 

As I promised, let me now in conclusion treat the general case, in which $\alpha$ is possibly less than $\theta$. 
In this case, I simply use instead the condition $p=\set{\<\kappa,\alpha>, \<\bar\kappa,\theta>}$ in the place of the earlier one, where $\bar\kappa$ is the next inaccessible beyond $\kappa$ and $\alpha$ in $N$. This condition ensures $j(f)(\kappa)=p(\kappa)=\alpha$ 
while still ensuring the closure of $\Ftail$. And since this was all that was required of $p$ in the previous argument, the theorem is proved.\QED

This Laver-like flexibility of the fast function makes it very useful in defining forcing iterations $\P_\kappa$, such as those in the next section. Our ability to specify $j(f)(\kappa)$ almost arbitrarily will allow us to exactly control the stage $\kappa$ forcing that appears in $j(\P)$ and the tail forcing after stage $\kappa$. It is this crucial control, the ability to push up the next nontrivial stage of forcing in $j(\P)$ after $\kappa$ arbitrarily high, which will allow us to develop a general lifting technique for unfoldability embeddings. 

\section{Unfoldability and the GCH}

I would like now to turn to the question of the failure of the GCH at an unfoldable cardinal. It is well-known that the consistency strength of the GCH failing at a measurable cardinal is strictly greater than the consistency strength of a measurable cardinal alone. So it is natural to wonder whether a similar fact is true for unfoldable cardinals. Villaveces \cite{Villaveces98} conjectured that the consistency strength of the failure of the GCH at an unfoldable cardinal is no greater than that of a strongly unfoldable cardinal, and he and Leshem \cite{Villaveces&Leshem} have announced a proof of this conjecture using class forcing over $L$. They report that Matt Foreman has also announced a proof that any strongly unfoldable cardinal $\kappa$ can be preserved to a forcing extension in which $2^\kappa>\kappa^\plus$, assuming a proper class of inaccessible cardinals above $\kappa$, by a proof, again, that uses class forcing. 

Here, I give a direct proof, over any model of ZFC, by set-forcing, that any unfoldable cardinal $\kappa$ can be made indestructible by the forcing to pump up $2^\kappa$ as large as desired. Since the method of proof seems quite general---I lift the unfoldability embeddings to the forcing extension---I expect that it will be adaptable to other sorts of forcing iterations with unfoldable cardinals.

\begin{Main Theorem} \label{Main} If $\kappa$ is unfoldable in $V$, then there is a forcing extension, by forcing of size $\kappa$, in which the unfoldability of $\kappa$ becomes indestructible by $\Add(\kappa,\theta)$ for any $\theta$.
\end{Main Theorem}

\Proof Villaveces has pointed out (Fact 3 in \cite{Villaveces98}) that if an unfoldable cardinal $\kappa$ is indestructible by $\Add(\kappa,1)$ then it is indestructible by $\Add(\kappa,\theta)$ for any $\theta$. So in order to prove the theorem it suffices merely to obtain a model in which the unfoldability of $\kappa$ is indestructible by $\Add(\kappa,1)$.

Suppose that $\kappa$ is unfoldable in $V$. Let $f\from\kappa\to\kappa$ be a $V$-generic fast function. The results of the previous section establish that $\kappa$ is unfoldable in $V[f]$, that every canonical $(\theta+1)$-unfoldability embedding in $V$ lifts to one in $V[f]$, and furthermore that we have complete freedom to specify the value of $j(f)(\kappa)$ for these lifted embeddings.

Define now in $V[f]$ the reverse Easton $\kappa$-iteration $\P$ which adds a Cohen subset to $\gamma$ at stage 
$\gamma$ for every $\gamma\in\dom(f)$, with trivial forcing at all other stages, and suppose that $G$ is $V[f]$-generic for $\P$. Let $A$ be a 
Cohen subset of $\kappa$ added generically over $V[f][G]$. The desired model will be $V[f][G][A]$. 

Fix any $\theta$ and any canonical $(\theta+1)$-unfoldability embedding $j:M\to N$ in the ground model $V$. I aim to lift this embedding to an embedding $j:M[f][G][A]\to N[j(f)][j(G)][j(A)]$ in the forcing extension. As an initial step, Theorem \ref{FastFunction} allows us to lift this embedding to $j:M[f]\to N[j(f)]$ with $j(f)(\kappa)>\theta$. The forcing $j(\P)$ factors as 
$\P*\Add(\kappa,1)*\Ptail$ and so we may use $G*A$ as the generic up to stage $\kappa$. The tail forcing $\Ptail$ is ${\leq}\theta$-closed in $N[j(f)][G][A]$, since the next nontrivial stage of forcing is the next element of $\dom(j(f))$, and this is neccessarily beyond $\theta$. As in the proof of Theorem \ref{FastFunction}, we factor: let $X=\set{j(g)(\kappa,\theta) \st g\in M[f] }$ and $\pi:X\iso N_0[j_0(f)]$, leading to the following diagram, where $j_0=\pi\compose j$ and $h=\pi^\inverse$.

\factordiagram{M[f]}{j_0}{j}{N_0[j_0(f)]}{h}{N[j(f)]}

Let me argue, in some detail, that $N_0[j_0(f)]^\ltkappa\of N_0[j_0(f)]$ in $V[f]$. Since every set in $N_0[j_0(f)]$, a model of $ZFC^-$, can be enumerated, a sequence of sets is equivalent in the model to a sequence of ordinals. So it suffices for me to show that $(\zeta^\ltkappa)^{V[f]}\of N_0[j_0(f)]$, where $\zeta=ORD^{N_0}=ORD^{N_0[j_0(f)]}$. Indeed, I will show that $(\zeta^\ltkappa)^{V[f]}\of N_0[f]$. Suppose that $s=\<\gamma_\alpha\st\alpha<\beta>\in(\zeta^\ltkappa)^{V[f]}$. Since $\beta<\kappa$, it follows that actually $s$ is in $V[f\restrict\beta]$, since the forcing to add the tail $f\restrict[\beta,\kappa)$ is $\leqbeta$-closed. Thus, $s$ has a $\F_\beta$-name $\dot s$ in $V$. This name gives rise to a sequence $\<\dot s_\alpha\st\alpha<\beta>$ of $\F_\beta$-names, where $\dot s_\alpha$ is a nice $\F_\beta$-name for the $\alpha^{th}$ ordinal in the sequence named by $\dot s$. Since the forcing $\F_\beta$ has size less than $\kappa$, a nice $\F_\beta$-name for a $\beta$-sequence or ordinals can itself be coded with a sequence of ordinals of length less than $\kappa$. Thus, since $N_0^\ltkappa\of N_0$ in $V$, each name $\dot s_\alpha$ is in $N_0$. Thus also the sequence $\<\dot s_\alpha\st\alpha<\beta>$ of names is also in $N_0$. And so, interpreting those names by the generic $f$, we conclude that $s=\<(\dot s_\alpha)_f\st\alpha<\beta>$ itself is in $N_0[f]$, as I had claimed. 

A similar argument establishes that $N_0[j_0(f)][G]^\ltkappa\of N_0[j_0(f)][G]$ in $V[f][G]$ and $N_0[j_0(f)][G][A]^\ltkappa\of N_0[j_0(f)][G][A]$ in $V[f][G][A]$. For these claims, it once again suffices to consider only $\zeta^\ltkappa$ in the appropriate model; and once again, a $\beta$ sequence of ordinals in $V[f][G][A]$ is actually in $V[f\restrict\beta][G\restrict\beta]$, and so the earlier argument shows that a name for a sequence of ordinals transforms to a sequence of names for ordinals in $N_0$. One then evaluates these names at the relevent generic object to obtain the original sequence in the desired model. 

Returning to the main line of argument, now, I may factor $j_0(\P)$ as $\P*\Add(\kappa,1)*\Potail$. Using $G*A$ as the generic up to stage $\kappa$, it remains to construct a tail generic $\Gotail\of\Potail$. This I will do by the same diagonalization technique as in Theorem \ref{FastFunction}. The first step is to observe that since $N_0$ has size $\kappa$, so also does $N_0[j_0(f)][G][A]$. Thus, in $V[f][G][A]$ I may line up the dense subsets of $\Potail$ in $N_0[j_0(f)][G][A]$ into a $\kappa$-sequence. I now construct recursively a descending $\kappa$-sequence of conditions which meet these dense sets one-by-one. At limit stages, the claim of the previous paragraph ensures that the construction up to that stage exists in $N_0[j_0(f)][G][A]$, and since $\Potail$ is ${<}\kappa$-closed in that model, there is a condition below it and the construction proceeds up to $\kappa$. The filter $\Gotail$ generated by this descending sequence of conditions exists in $V[f][G][A]$, but because it meets the relevent dense sets, it is $N_0[j_0(f)][G][A]$-generic for $\Potail$.

The embedding $j_0$ therefore lifts to $j_0:M[f][G]\to N_0[j_0(f)][j_0(G)]$ with $j_0(G)=G*A*\Gotail$. Since the 
critical point of $h$ is bigger than $\kappa$, the embedding $h$ lifts trivially to 
$h:N_0[j_0(f)][G][A]\to N[j(f)][G][A]$. To lift $h$ through the rest of the $j_0(\P)$ forcing, I will employ an argument similar the previous case to show that the filter $\Gtail$ generated by $h\image\Gotail$ is $N[f][G][A]$-generic for $\Ptail$. Namely, any dense set $D$ has 
the form $h(r)(s)$ for some $r:\theta_0^{|s|}\to N_0[j_0(f)][G][A]$ with $r\in N_0[j_0(f)][G][A]$ and some $s\in\theta^{{<}\omega}$. I may assume without loss of generality that $r(\sigma)$ is dense in $\Potail$ for every $\sigma$. Let $D^*$ be the intersection of all $r(\sigma)$. Thus, $h(D^*)$ is contained in $D$. But by 
the genericity of $\Gotail$, we know $D^*$ meets $\Gotail$, and so $h(D*)$ 
meets $h\image \Gotail$. Thus, $D$ meets $h\image\Gotail$, as desired. Putting everything together, I conclude that $G*A*\Gtail$ is $N[j(f)]$-generic for $j(\P)$. So we 
can lift $j$ to $j:M[f][G]\to N[j(f)][j(G)]$ with $j(G)=G*A*\Gtail$, resulting in the following diagram:

\factordiagram{M[f][G]}{j_0}{j}{N_0[j_0(f)][j_0(G)]}{h}{N[j(f)][j(G)]}

It remains to lift through the $A$ forcing. Here, we observe that $A$, being a bounded subset of $j(\kappa)$ in $N[j(f)][j(G)]$, is a 
condition in $j(\Add(\kappa,1))$. The forcing $j(\Add(\kappa,1))$ is $\leqtheta$-closed, so we may use the same 
factor argument to diagonalize to produce an 
$N_0[j_0(f)][j_0(G)]$-generic $A_0\of j_0(\kappa)$ below the condition $A$. Again, 
$h\image A_0$ will 
generate an $N[j(f)][j(G)]$-generic filter $j(A)$ for $j(\Add(\kappa,1))$. So we may 
lift the embedding to $j:M[f][G][A]\to N[j(f)][j(G)][j(A)]$, as desired. 

\factordiagram{M[f][G][A]}{j_0}{j}{N_0[j_0(f)][j_0(G)][j_0(A)]}{h}{N[j(f)][j(G)][j(A)]}

At this point, I have established that any canonical $(\theta+1)$-unfoldability embedding $j:M\to N$ in the ground model $V$ can be lifted to a $(\theta+1)$-unfoldability embedding $j:M[f][G][A]\to N[j(f)][j(G)][j(A)]$ in the forcing extension. It follows that $\kappa$ is unfoldable in $V[f][G][A]$, because if $M'$ is a structure of size $\kappa$ in $V[f][G][A]$, then it has a name $\dot M'$ of size $\kappa$ in $V$, and for any $\theta$ this name can be placed into such an $M$ for which I will have lifted the embedding. The structure $M'=(\dot M')_{f*G*A}$ therefore exists in $M[f][G][A]$, and the restriction $j\restrict M'$ provides an unfoldability witness embedding for $M'$. 

So $\kappa$ is unfoldable in $V[f][G][A]$. The final step is to realize that since adding two Cohen subsets to $\kappa$ is the same as adding one, if we force over $V[f][G][A]$ to add another Cohen subset $A'\of\kappa$, then we may regard $A*A'$ as a single Cohen subset of $\kappa$ and employ the previous argument to conclude that $\kappa$ is unfoldable in $V[f][G][A*A']$. So the unfoldability of $\kappa$ is indestructible in $V[f][G][A]$ by $\Add(\kappa,1)$, as desired.\QED

A minor modification to the argument allows one to use $V[f][G]$ as the desired model, rather than $V[f][G][A]$, and avoid the argument involving $A'$ at the conclusion of the proof. Specifically, to accomplish this one should modify $\P$ so that forcing is only done at stages $\gamma\in\dom(f)$ for which $f(\gamma)$ is an successor ordinal. By ensuring that $j(f)(\kappa)$ is a sufficiently large limit ordinal, there will be trivial forcing at stage $\kappa$ in $j(\P)$, and one may conclude that $\kappa$ is unfoldable in $V[f][G]$. Subsequently, in $V[f][G][A]$, one takes $j(f)(\kappa)$ to be a sufficiently large successor ordinal and follows the previous argument. With this modification, consequently, $\kappa$ will be unfoldable in $V[f][G]$ and indestructible from there to $V[f][G][A]$. 

As I promised earlier, let me explain how this proof fits into the context of Theorem 2.1 of \cite{Villaveces&Leshem}. Supposing $V=L$, we began with a canonical $(\theta+1)$-unfoldability embedding $j:M\to N$. We then lifted this embedding to 
$$j:M[f][G][A]\to N[j(f)][j(G)][j(A)].$$
For any $S\of\kappa$ in $M[f][G][A]$, therefore, we have $$\left\langle\strut(L[f][G][A])_\kappa,{\in},S\right\rangle\elesub\left\langle\strut(N[j(f)][j(G)][j(A)])_{j(\kappa)},{\in},j(S)\right\rangle,$$
witnessing the $\theta$-unfoldability of $\kappa$ in $L[f][G][A]$ according to Villaveces' original definition. Villaveces and Leshem, in the course of their proof of Theorem 2.1 in \cite{Villaveces&Leshem}, claim, in their terminology, that $M\of L_\rho[G]$; here, that claim translates to the assertion that $L_\rho[j(f)][j(G)][j(A)]\of L_\rho[f][G][A]$, where $\rho=j(\kappa)$ is the height of the extending structure $(N[j(f)][j(G)][j(A)])_{j(\kappa)}$. But this claim cannot be correct, since $f$, $G$ and $A$ are already contained in $L_\rho[j(f)][j(G)_{\kappa+1}]$, and $j(G)$ adds further generic objects over this model. Indeed, this same feature arises already with the embedding $j:M[f]\to N[j(f)]$, where we have $L_\kappa[f]\elesub L_{j(\kappa)}[j(f)]$ and $L_{j(\kappa)}[j(f)]\of L[f]$, but $j(f)$ adds objects above $\kappa$ that are generic over over $L_{j(\kappa)}[f]$, so $L_{j(\kappa)}[j(f)]\not\of L_{j(\kappa)}[f]$. The point is that one needs ordinals larger than $\rho=j(\kappa)$ in order to define the embeddings $h$ and $j$ used in the construction of the generic objects $j(f)$, $j(G)$ and $j(A)$. In summary, I believe that the proof of  \cite{Villaveces&Leshem} Theorem 2.1 is incorrect on this point and the statement of their theorem is proved false by the Main Theorem of this paper.

\Theorem. If $\kappa$ is strongly unfoldable in $V$, then this is preserved by the forcing in Theorems \ref{FastFunction} and \ref{Main}. In particular, any strongly unfoldable cardinal $\kappa$ can be made indestructible by $\Add(\kappa,1)$.\label{Strong}

\Proof This is a simple, general observation. If $\kappa$ is strongly unfoldable in $V$, then for every $\theta$ there is a canonical $\beth_\theta+1$-unfoldability embedding $j:M\to N$ with $V_\theta\of N$. In the case of Theorem \ref{FastFunction}, this embedding lifts to $j:M[f]\to N[j(f)]$. Since $f\in N[j(f)]$, it follows that $V_\theta[f]\of N[j(f)]$, and so the lifted embeddings witness the strong unfoldability of $\kappa$ in $V[f]$. Similarly, in the case of the Main Theorem \ref{Main}, the embedding lifts to $j:M[f][G][A]\to N[j(f)][j(G)][j(A)]$, and since, once again, both $V_\theta$ and the generics $f*G*A$ are in $N[j(f)][j(G)][j(A)]$, the lifted embedding remains sufficiently strong. The argument involving $A'$ then shows that further forcing with $\Add(\kappa,1)$ must preserve the strong unfoldability of $\kappa$.\QED

While the previous theorem shows that any strongly unfoldable cardinal $\kappa$ can be made indestructible by $\Add(\kappa,1)$, it does not immediately follow from this, as it does for mere unfoldability, that the strong unfoldability of $\kappa$ is also indestructible by $\Add(\kappa,\theta)$ for any $\theta$. One might ask the question:

\Question. If the strong unfoldability of $\kappa$ is indestructible by $\Add(\kappa,1)$, is it also indestructible by $\Add(\kappa,\theta)$ for every $\theta$?

This is answered, in the negative, by the following theorem. 

\Theorem. If $\kappa$ is strongly unfoldable, then there is a forcing extension, by forcing of size $\kappa$, in which the strong unfoldability of $\kappa$ is indestructible by $\Add(\kappa,1)$ but not by $\Add(\kappa,\theta)$ for any cardinal $\theta>\kappa$.

\Proof Use the model $V[f][G][A]$ of the Main Theorem. I have proved already that if $\kappa$ is strongly unfoldable in $V$, then it remains so in $V[f][G][A]$. It follows by the argument involving $A'$ at the conclusion of the proof of the Main Theorem, that the strong unfoldability of $\kappa$ is indestructible by further forcing with $\Add(\kappa,1)$. Now suppose that $\vec A$ is $V[f][G][A]$-generic for further forcing by $\Add(\kappa,\theta)$ for some cardinal $\theta>\kappa$, and, towards a contradiction, that $\kappa$ remains strongly unfoldable in $V[f][G][A][\vec A]$. Suppose that $M$ is a nice model of $ZFC^-$ in $V$ of size $\kappa$. It follows that $M[f][G][A]$ is a nice model of $ZFC^-$ in $V[f][G][A]$. Suppose that $j:M[f][G][A]\to\tilde N$ is a strong $\kappa+1$-unfoldability embedding. Neccessarily, $\tilde N$ has the form $N[j(f)][j(G)][j(A)]$. 

I claim that $P(\kappa)^{N[j(f)]}\of V[f]$. First, I observe that $N_\kappa=V_\kappa$, and so since the tail forcing is closed, $(N[j(f)])_\kappa=(V[f])_\kappa$. Thus, if $A\of\kappa$ and $A\in N[j(f)]$, then every initial segment of $A$ is in $(N[j(f)])_\kappa$, and hence in $V[f]$. Since the iteration $\P$ admits a gap below $\kappa$---that is, it factors as $\P_1*\P_2$ where $\card{\P_1}<\gamma$ and $\forces\P_2$ is $\leqgamma$-closed for some $\gamma<\kappa$---it follows by the Key Lemma of \cite{GapForcing} that it cannot add any new subsets of $\kappa$ all of whose initial segments are in the ground model. So $A$ must be in $V[f]$. 

It follows from this, in particular, that none of the sets in $\vec A$ are in $N[j(f)]$. Furthermore, since the forcing $j(\P)_{\kappa+1}$ has size $\kappa$, it can add at most $\kappa$ many of the sets in $\vec A$ (and these are added at stage $\kappa$ if at all). Since the rest of the forcing $j(\P)_{\kappa+1,j(\kappa)}*\Add(j(\kappa),1)$ does not add any additional subsets to $\kappa$, most of the sets in $\vec A$ are not in $N[j(f)][j(G)][j(A)]$. This contradicts the fact that the embedding was supposed to be $\kappa+1$-strongly unfoldable.\QED

\Corollary. By forcing of size $\kappa^\plus$, one can make any unfoldable cardinal $\kappa$ unfoldable but not strongly unfoldable.

\Proof Simply force to $V[f][G][A]$, and then force with $\Add(\kappa,\kappa^\plus)$. The previous arguments establish that this forcing will preserve the unfoldability but not the strong unfoldability of $\kappa$.\QED

\section{Universal indestructibility}\label{Universal}

In this last section, I would like to treat all unfoldable cardinals simultaneously, making the GCH fail at all of them. Villaveces and Leshem \cite{Villaveces&Leshem} have proved, by forcing over $L$, one can make the GCH fail at every inaccessible cardinal while preserving the strong unfoldability of every unfoldable cardinal in $L$. Here, I prove that over any model of ZFC one can force the GCH to fail at every inaccessible cardinal (or, just as well, every regular cardinal) while preserving every strongly unfoldable cardinal. Since in $L$ every unfoldable cardinal is strongly unfoldable, this theorem therefore generalizes, by a different proof, the \cite{Villaveces&Leshem} result. Using the embedding characterization of unfoldability, I will simply lift the strong unfoldability embeddings to the forcing extension.

Before proving the theorem, let me remark that it is not difficult to prove that $\kappa$ is strongly unfoldable if and only if for every set $S\of\kappa$ and every $\theta$ there is a nice model $M$ with $S\in M$ and an embedding $j:M\to N$ with $V_\theta\of N$ and $j(\kappa)\geq\theta$. The reason, using the definition of \cite{Villaveces&Leshem}, is that from this it follows that $\<V_\kappa,{\in},S>\elesub \<N_{j(\kappa)},{\in},j(S)>$ and $V_\theta\of N_{j(\kappa)}$. 

\Theorem. There is a forcing extension preserving all strongly unfoldabile cardinals, in which the GCH fails at every inaccessible cardinal.

\Proof Suppose that $G$ is $V$-generic for the (possibly proper class) reverse Easton iteration $\P$ that forces with $Q_\gamma=\Add(\gamma,\gamma^\plusplus)$ at every inaccessible stage $\gamma$, and trivial forcing at all other stages. It has been argued elsewhere that $V[G]$ is a model of ZFC. I will show that every strongly unfoldable cardinal $\kappa$ from $V$ is preserved to $V[G]$. In order to do so, I will lift the canonical ground model strong unfoldability embeddings $j:M\to N$ to strong unfoldability embeddings $j:M[G_\kappa][A]\to N[j(G_\kappa)][j(A)]$ in the extension, where $A\of\kappa$ is any Cohen subset of $\kappa$ in $V[G]$. This suffices because any set $S\of\kappa$ in $V[G]$ is actually in $V[G_\kappa][A]$ for some such $A$, and therefore in $M[G_\kappa][A]$ for a suitably large choice of $M$. 

Choose any non-cardinal $\theta\muchgt\kappa$ and any canonical strong unfoldability embedding $j:M\to N$ in the ground model, where $M$ is a nice model of $ZFC^-$, $V_\theta\of N$, $j(\kappa)>\theta$ and $N=\set{j(g)(s)\st g\in M\And s\in\delta^\ltomega}$ where $\delta=\beth_\theta=\card{V_\theta}$. Consider the forcing $j(\P)$ in $N$. Since $V_\theta\of N$, the iteration up to stage $\theta$ in $N$ and $V$ is the same, and I may factor $j(\P)$ as $\P_\theta*\Ptail$, where $\Ptail$ is the forcing beginning from stage $\theta$. And since $G_\theta$ is $V$-generic for $\P_\theta$, it is also $N$-generic. Thus, in order to lift the embedding, I must merely construct in $V[G]$ an $N[G_\theta]$-generic for $\Ptail$. Since $\theta$ is not a cardinal, the first nontrivial stage of forcing in $\Ptail$ occurs at the next innaccessible of $N$ beyond $\theta$, and this is larger than $\delta=\beth_\theta=(\beth_\theta)^N$. Thus, $\Ptail$ is $\leqdelta$-closed in $N[G_\theta]$. 

Let $X=\set{j(g)(\kappa,\theta)\st g\in M}$. My earlier arguments establish that $X\elesub N$ and $X^\ltkappa\of X$ in $V$. Furthermore, since $M$ has size $\kappa$, there are only $\kappa$ many functions $g$ to represent elements of $X$, and so $X$ has size $\kappa$. Now let $Y=X[G_\theta]=\set{\tau_{G_\theta}\st \tau\in X}$ be the set of interpretations of $\P_\theta$-names in $X$ by the generic object $G_\theta$. (Note: since $G_\theta$ may not be $X$-generic, it may be that $Y$ has more ordinals than $X$.) The set $Y$ also has size $\kappa$. 

I claim that $Y\elesub N[G_\theta]$. To see this, simply verify the Tarski-Vaught criterion: if $N[G_\theta]\satisfies\exists x\,\varphi(\tau_{G_\theta},x)$ for some $\tau\in X$, then the boolean value $b=\boolval{\exists x\,\varphi(\tau,x)}$, which is in $X$, is non-zero (and met by $G$). Thus, there is some name $\sigma\in X$ with $b=\boolval{\varphi(\tau,\sigma)}$. Since $G$ meets $b$, it follows that $N[G]\satisfies\varphi(\tau_{G_\theta},\sigma_{G_\theta})$, and the criterion is verified. 

Next, I claim that $Y^\ltkappa\of Y$ in $V[G_\theta]$. Since any set can be enumerated, it suffices to show that $(Y\intersect ORD)^\ltkappa\of Y$ in $V[G_\theta]$. So suppose that $s\in V[G_\theta]$ is a $\beta$-sequence of ordinals in $Y$, with $\beta<\kappa$. By the closure of the tail forcing, we know that actually $s\in V[G_{\beta+1}]$, and so $s$ has a $\P_{\beta+1}$-name $\dot s\in V$. This name gives rise in $V$ to a sequence $\<\dot s_\alpha\st\alpha<\beta>$ of nice $\P_{\beta+1}$-names for ordinals. Since such names can be coded with a $\beta$-sequence of ordinals, they are all in $X$. Consequently, the whole $\beta$-sequence of these names must also be in $X$. Thus, the sequence of interpretations of these names by the generic $G_\theta$, namely the sequence $s$, must be in $Y=X[G_\theta]$, as desired. 

Enumerate in $V[G]$ the $\kappa$ many dense subsets of $\Ptail$ in $Y$, and construct a descending $\kappa$-sequence of conditions that meet these dense sets one-by-one. At limit stages, the sequence constructed so far lies in $Y$ because $Y^\ltkappa\of Y$, and since the forcing $\Ptail$ is $\leqtheta$-closed, the construction may continue up to $\kappa$. Let $\Gtail$ be the filter generated by the resulting descending sequence. By construction, $\Gtail$ is $Y$-generic. 

I claim that $\Gtail$ is neccessarily $N[G_\theta]$-generic. To see this, suppose $D\of\Ptail$ is an open dense subset of $\Ptail$ in $N[G_\theta]$. The set $D$ must have the form $j(\vec D)(s)_{G_\theta}$ for some sequence $\vec D\in M$ of names and some other parameter $s\in\delta^\ltomega$. Let $\bar D$ be the intersection of all $j(\vec D)(t)_{G_\theta}$ which are open dense subsets of $\Ptail$ for any $t\in \delta^\ltomega$. Thus $\bar D\of D$ and since $\Ptail$ is $\leqdelta$-closed in $N[G_\theta]$, the set $\bar D$ remains dense. Furthermore, since by intersecting over all $t$ we have eliminated the need for the parameter $s$, it follows that $\bar D\in Y$. Consequently, since $\Gtail$ is $Y$-generic, it meets $\bar D$, and hence also $D$, as desired. 

I may therefore lift the embedding to $j:M[G_\kappa]\to N[j(G_\kappa)]$ where $j(G_\kappa)=G_\theta*\Gtail$. Suppose now that $A\of\kappa$ is one of the generic Cohen sets added at stage $\kappa$ by $G$. By rearranging the sets in $G$ if neccessary, I may assume that $A$ is the first set added at stage $\kappa$ in $G$, and hence also in $j(G_\kappa)$. In particular, $A$ is a bounded subset of $j(\kappa)$ in $N[j(G_\kappa)]$, and therefore a condition in $j(\Add(\kappa,1))$ in that model. By the same technique just used to construct $\Gtail$, I may construct a filter below the condition $A$ in $j(\Add(\kappa,1))$ that is generic over $Y'=X[j(G_\kappa)]$. And once again, the argument shows that this filter will be fully $N[j(G_\kappa)]$-generic. Thus, I may lift the embedding to $j:M[G_\kappa][A]\to N[j(G_\kappa)][j(A)]$, where $j(A)$ is the (union of) the generic filter I just mentioned. Since $V_\theta\of N$ and $G_\theta\in N[j(G_\kappa)]$, it follows that $(V[G])_\theta\of N[j(G_\kappa)][j(A)]$, so this lifted embedding is a $\theta$-strong unfoldability embedding. 

Since $\theta$ can be taken to be arbitrarily large, it follows that $\kappa$ is strongly unfoldable in $V[G]$. And since $\P$ clearly forces the failure of the GCH at every inaccessible cardinal, the theorem is proved.\QED

The proof is highly mallable. In particular, it is not neccessary to restrict the forcing to the inaccessible cardinal stages. Indeed, while preserving every unfoldable cardinal, one can force the GCH to fail at every regular cardinal $\gamma$. What is more, it is an easy matter to modify the proof to ensure $2^\gamma=\gamma^{+++}$ at all these cardinals, rather than simply $2^\gamma=\gamma^\plusplus$, or whatever is desired. Specifically, for any suitably definable function $f:ORD\to ORD$, one may force $2^\gamma=\aleph_{f(\gamma)}$ for every regular cardinal $\gamma$ while preserving every strongly unfoldable cardinal $\gamma$. In addition to the usual Easton requirements on $f$, one need only require that the unfoldable cardinals are closed under $f$ and that the function is defined locally enough so that if $V_\theta\of N$, then $j(f)$ in $N$ agrees with $f$ up to $\theta$, so that the forcing in $j(\P)$ in $N$ up to stage $\theta$ is the same as that in $\P$. 

In closing this paper, let me remark on a curious point I have stumbled in the course of proving the theorems in this paper, namely, the question of whether of the Levy-Solovay theorem holds for unfoldable cardinals. Of course, anyone would expect that it does, but a proof is not forthcoming. Specifically, the question is whether the unfoldability of $\kappa$ in a forcing extension by forcing of size less than $\kappa$ implies the unfoldability of $\kappa$ in the ground model. One could generalize this question in light of my recent work on gap forcing (see \cite{GapForcingGen}, \cite{GapForcing}), and ask if the same is true in any forcing extension admitting a gap below $\kappa$. That is, if $\kappa$ is unfoldable in $V^{\P*\Qdot}$ where $\card{\P}=\gamma<\kappa$ and $\forces\Qdot$ is $\leqgamma$-strategically closed, must $\kappa$ be unfoldabe in $V$? 

\bibliography{Unfold}
\bibliographystyle{alpha}
\end{document}